\documentclass{article}
\usepackage[dvips]{graphicx}
\usepackage{amssymb,amsfonts,amsmath}
\usepackage{amsmath,amsmath,amsfonts}
\usepackage{amssymb}
\usepackage{epsfig}
\usepackage[figureright]{rotating}

\begin{document}

\title{Interpolation and Iteration for Nonlinear Filters}
\author{Alexandre J.\ Chorin and Xuemin Tu\\\\
Department of Mathematics, \\
University of California at Berkeley
and
Lawrence Berkeley National Laboratory,\\
Berkeley, CA, 94720}

\date{}
\maketitle

\begin{abstract}
We present a general form of the iteration and interpolation process used
in implicit particle filters. Implicit filters are based on a pseudo-Gaussian representation of
posterior densities, and are designed to focus the particle paths so as to reduce the
number of particles needed 
in nonlinear data assimilation. 
Examples are given. 
\end{abstract}

%
%
{\bf Keywords:} {Implicit sampling, filter, pseudo-Gaussian, Jacobian, chainless sampling, particles}

\section{Introduction}

There are many problems in science in which the state
of a system must be identified from an uncertain equation supplemented by a stream of noisy data (see e.g. \cite{Do2}). A natural model of this situation consists of an Ito stochastic differential equation (SDE):
\begin{equation}
dx=f(x,t) \, dt+  g(x,t) \, dw,
\label{eq:datass}
\end{equation}
where $x=(x_1,x_2,\dots,x_m)$ is an $m$-dimensional vector, $w$ is $m$-dimensional Brownian motion, $f$ is an $m$-dimensional vector function, and 
$ g(x,t)$ is an $m$ by $m$ diagonal matrix.  
The initial state
$x^0$ is assumed given and may be random as well.

As the solution of the SDE unfolds, it is observed, and the
values $b^n$ of a measurement process are recorded at times $t^n, n=1,2,...$
For simplicity assume $t^n=n\delta$, where $\delta$ is a fixed
time interval.
The measurements are related to the evolving state $x(t)$ by
\begin{equation}
b^n={h}(x^n)+Q W^n,
\label{eq:observe}
\end{equation}
where $h$ is a $k$-dimensional, generally nonlinear, vector function with $k \le m$, $Q$ is a $k$ by $k$ diagonal matrix, $x^n= x(n\delta)$,  and $ W^n$ is a vector whose
components are $k$ independent Gaussian variables of mean zero and variance one, independent also of the Brownian motion in equation (\ref{eq:datass}). The task is to estimate $x$ on the basis of equation (\ref{eq:datass}) and the observations (\ref{eq:observe}).

If the system~(\ref{eq:datass}) and equation (\ref{eq:observe}) are linear and the data are Gaussian, the solution can be found via the Kalman-Bucy filter (see e.g. \cite{BSM}). In the general case, it is natural to try to estimate $x$ via its evolving probability density.
The initial state $x^0$ is known and so is its probability density; all one has to do is
evaluate sequentially the density $P_{n+1}$ of $x^{n+1}$ given the probability
densities $P_k$ of
$x^k$ for $k\le n$ and the data $b^{n+1}$. This can be done by following
``particles" (replicas of the system) whose empirical distribution approximates $P_n$. 
A standard construction 
(see e.g 
\cite{Mac1,LS,Dou1,Ar2,Gi1,Ch11,Dow1,toapp})
uses the probability density function (pdf) $P_n$ and equation (\ref{eq:datass}) to generate
a prior density, and then uses the new data $b^{n+1}$ to generate a
posterior density $P_{n+1}$ through weighting and resampling. In addition, one has to sample backward to take into account the information
each measurement provides about the past, as well as avoid having too many identical particles after resampling. This can be very expensive, in particular because 
the number of particles needed can grow catastrophically
(see e.g. 
\cite{Sny,Blb} and also Example 2 below). Sophisticated methods for generating efficient priors
can be found e.g. in \cite{Dou1,Ar2}. The challenge is to generate high
probability samples 
so as to minimize the effort of computing particle paths whose weight is very low. 

In \cite{CT1} we introduced an alternative to the standard approach. In our 
method the posterior density is
sampled directly by iteration and interpolation, as suggested by our earlier work on chainless sampling \cite{Ch101}, and by the observation in \cite{We1} connecting interpolation and the marginalization process used in chainless sampling. The new filter aims the particle
 trajectories 
as accurately as possible in the direction of the observations so that fewer particles are needed. In that earlier paper our approach was presented by means of simple examples.
In the present paper we present a general, more abstract, formulation, introduce an extension to the
case of sparse observations, and discuss additional examples.

\section{Forward step}\label{sec_forward}
To begin, assume
that at time $t^n=n\delta$, where $\delta>0$ is fixed, we have
a collection of $M$ particles $X_i^n$, $1 \le i \le M$, $n=0,1,\dots$, whose
empirical density approximates $P_n$, the probability density at time $n\delta$
of the particles that obey the evolution equation (\ref{eq:datass}) subject to the
observations (\ref{eq:observe}) at times $t=k\delta$ for $k \le n$. In the present section we explain how to
find
positions for the same particles at time $(n+1)\delta$ given only the positions at time $n\delta$ and the pdf $P_n$, taking into
account the next observation and the equation of motion. 
Let $N(a,v)$ denote a Gaussian variable of mean $a$ and variance $v$.
First, approximate the SDE (\ref{eq:datass}) by a difference scheme of the form
\begin{equation}
X^{n+1}=X^n+F(X^n,t^n)\delta+G(X^n,t^n)V^{n+1},
\label{appeq}
\end{equation}
where we assume temporarily that $\delta$ equals the interval between observations, 
i.e., we assume that there is an observation at every time step.
$X^n$ stands for $X(n\delta)$, $G$ is assumed to be diagonal, and $X^n,X^{n+1}$
are $m$ dimensional vectors. 
$F,G$ determine the scheme used to solve the SDE, see for example \cite{CT1}. 
$V^{n+1}$ is a vector of $N(0,\delta)$ Gaussian variables, independent of 
each other for each $n$, with the vectors $V^{n+1}$ independent of each other for
differing $n$, independent also of the $W^k, k=1,...,$ in the observation equation 
(\ref{eq:observe}).
The sequence of $X^n, n=0,1,\dots$ approximates a sample solution of the SDE,
$X^0$ is assumed given and may be random. 
The function $G$ in (\ref{appeq}) does not depend on $X^{n+1}$ for an Ito equation,
and we assume for simplicity that $F$ does not depend on $X^{n+1}$ either,
because this was the case in all the examples we have worked on so far. The analysis
below can be easily repeated for the case where $F$ does depend on $X^{n+1}$,
at the cost of slightly more complicated formulas. 
Equation (\ref{appeq}) states that $X^{n+1}-X^n$ is 
an $N(F(X^n,t^n)\delta,\delta G(X^n,t^n)^*G(X^n,t^n))$ vector, where the star * denotes a transpose.

We have one sample solution $X_i^n$ of the SDE for each particle.
Our task is to sample, for each particle, the vector 
$X^{n+1}_i$ whose probability density is determined
by the approximation of the SDE as well as by the next observation 
for each of the $M$ particles. 
We keep the
notation $X^{n+1}_i$  for the positions of the particles even 
though once the observation is taken into account these positions 
no longer coincide with the positions
of sample solutions of equation (\ref{appeq}).

Consider the $i$-th particle. We are going to work particle by particle, so that the particle 
index $i$ will be temporarily suppressed.
Suppose we already know the posterior vector $X^{n+1}$. Its
probability density $P_{n+1}$ of $X^{n+1}$ given $X^n$ is
\begin{eqnarray}
P_{n+1}(X^{n+1})&=&Z^{-1}\exp\left(-\left(X^{n+1}-X^n-F_n\right)^*(G_n^*G_n)^{-1}\left(X^{n+1}-X^n-F_n\right)/2\right.\nonumber\\
&&\left.-\left(h(X^{n+1})-b^{n+1}\right)^*(Q^*Q)^{-1}\left(h(X^{n+1})-b^{n+1}\right)/2\right),\nonumber\\
\label{Pn+1}
\end{eqnarray}
where the functions $F_n=F(X^n,t^n)\delta$, and $G_n=\sqrt{\delta}G(X^n,t^n)$ can be read from the approximation of the SDE, and
$Z$ is a normalization constant, the integral of the numerator over
all $X^{n+1}$ with $X^n$ fixed. The value of this $Z$ is not 
available.
Our goal is to find samples $X^{n+1}$ whose probability is high, and which 
are well
distributed with respect to $P_{n+1}$. We do that by picking the 
probability in
advance: we first pick samples of $m$ $N(0,1)$ variables 
$(\xi_1,\xi_2,\dots,\xi_m)=\xi$,
whose joint pdf (probability density function) is
$\exp(-\xi^*\xi/2))/(2\pi)^{m/2}$, 
and require that each $X^{n+1}$ be a function of a sample $\xi$ with
 the same probability as $\xi$,
up to the Jacobian of the transformation. This should produce 
likely and well-distributed samples.

A little thought shows that this can be done,
not by equating $P_{n+1}$ to  $\exp(-\xi^*\xi/2)/(2\pi)^{m/2}$, but 
by equating the arguments of the two exponentials. For example,
if one wants to represent a $N(0,v)$ random variable $x$ with pdf 
$\exp(-\frac{x^2}{2v})/\sqrt{2\pi v}$
as a function of a $N(0,1)$ variable $\xi$ with pdf $\exp(-\xi^2/2)/\sqrt{2\pi}$, equating the
arguments yields $x=\sqrt{v}\,\xi$, clearly a good choice. 
Thus, we wish to solve the equation
\begin{flalign}
&{\xi}^*\xi/2=\nonumber\\
=&
\left(X^{n+1}-X^n-F_n\right)^*(G_n^*G_n)^{-1}\left(X^{n+1}-X^n-F_n\right)/2
+\left(h(X^{n+1})-b^{n+1}\right)^*(Q^*Q)^{-1}\left(h(X^{n+1})-b^{n+1}\right)/2\nonumber\\
\label{args}
\end{flalign}
and obtain $X^{n+1}$ as a function of $\xi$.

We proceed point by point--- given a vector $\xi$, we find the corresponding $X^{n+1}$ rather than
look for an expression for the function $X^{n+1}(\xi)$ as a whole---and by iteration: we find 
a sequence of approximations $X^{n+1}_j$ ($=X_j$ for brevity) which converges to $X^{n+1}$;
we set $X_0=0$, and now explain how to find $X_{j+1}$ given $X_j$.
First, expand the function $h$ in the observation equation (\ref{eq:observe}) in Taylor series around $X_j$:
\begin{equation}\label{lineeqn}
h(X_{j+1})=h(X_j)+H_{j}\cdot(X_{j+1}-X_j),
\end{equation}
where $H_{j}$ is a Jacobian matrix evaluated at $X_j$.
The observation equation (\ref{eq:observe}) can be approximated as:
\begin{equation}
z_j=H_{j}X_{j+1}+QW^{n+1},
\label{lincon}
\end{equation}
where $z_j=b^{n+1}-h(X_j)+H_{j}X_j$.

The left side of equation (\ref{args}) can be approximated as:
\begin{flalign}
&\left(X_{j+1}-X^n-F_n\right)^*(G_n^*G_n)^{-1}\left(X_{j+1}-X^n-F_n\right)/2
+\left(H_{j}X_{j+1}-z_j\right)^*(Q^*Q)^{-1}
\left(H_{j}X_{j+1}-z_j\right)/2\nonumber\\
=&\left(X_{j+1}-\bar{m}_j\right)^*\Sigma_j^{-1}\left(X_{j+1}-\bar{m}_j\right)/2+\Phi_j,
\label{H12}
\end{flalign}
where 
$$
\Sigma_j^{-1}=(G_n^*G_n)^{-1}+H_{j}^*(Q^*Q)^{-1}H_{j},\quad
\bar{m}_j=\Sigma_j\left((G_n^*G_n)^{-1}(X^n+F_n)+H_j^*(Q^*Q)^{-1}z_j\right),
$$
and
$$K_j=H_{j}G_n^*G_nH_{j}^*+Q^*Q,\quad
\Phi_j=\left(z_j-H_{j}(X^n+F_n)\right)^*K_j^{-1}\left(z_{j}-H_{j}(X^n+F_n)\right)/2.
$$

We now solve for $X_{j+1}$ as a function of $\xi$. To make the computation tractable, in this step we ignore the remainder $\Phi_j$;
this is a key step.  
We thus solve
the simpler equation
\begin{equation}
(X_{j+1}-\bar{m}_j)^*\Sigma_j^{-1}(X_{j+1}-\bar{m}_j)/2=\xi^*\xi/2.
\label{simpler}
\end{equation}
This can be done in any of a number of ways; for example, one can write $\Sigma_j=L_jL_j^*$, where
$L_j$ is a lower triangular matrix and $L_j^*$ is its transpose,
and then set $X_{j+1}=\bar{m}_j+L_j\xi$ (a different algorithm was suggested in \cite{CT1}).
The iteration is done. 

If the sequence $X_j$ converges to a limit, call the limit $X^{n+1}$. One can 
readily check that the approximate equation (\ref{lincon})
converges to the full observation equation (\ref{eq:observe}). The remainders
$\Phi_j$ also converge to a limit $\Phi^{n+1}$. Equation
(\ref{args}) becomes:
\begin{flalign}
&\xi^*\xi/2+\Phi^{n+1}=\nonumber\\
=&\left(X^{n+1}-X^n-F_n\right)^*(G_n^*G_n)^{-1}\left(X^{n+1}-X^n-F_n\right)/2
+(h(X^{n+1})-b^{n+1})(Q^*Q)^{-1}(h(X^{n+1})-b^{n+1})/2\nonumber.\\
\label{wemiss}
\end{flalign}
Multiply this equation by $-1$ and exponentiate both sides:
\begin{flalign}
&\exp(-\xi^*\xi/2)\exp(-\Phi^{n+1})=\nonumber\\
=&\exp\left(-\left(X^{n+1}-X^n-F_n\right)^*(G_n^*G_n)^{-1}\left(X^{n+1}-X^n-F_n\right)/2
-\left(h(X^{n+1})-b^{n+1})^*(Q^*Q)^{-1}(h(X^{n+1})-b^{n+1}\right)/2\right).\nonumber\\
\label{main}
\end{flalign}
This differs from what we set out to do in equation (\ref{args}) by the factor $\exp(-\Phi^{n+1})$ on the right hand side. 

Let $P(\alpha|\beta)$ be the probability of $\alpha$ given $\beta$. The factor 
$\exp(-\Phi^{n+1})$ is proportional to $P(b^{n+1}|X^n)$, and 
equation (\ref{main}) is the statement
\begin{equation}
P(X^{n+1}|X^n,b^{n+1})P(b^{n+1}|X^n)=P(X^{n+1}|X^n)P(b^{n+1}|X^{n+1}),
\end{equation}
i.e., this is Bayes' theorem. 
Note also that 
 equation (\ref{simpler}) is a pseudo-Gaussian representation of
$X^{n+1}$, not a Gaussian representation; the matrix $\Sigma_j$ is a function of the sample. 

We next compute the Jacobian determinant $J=\det({\partial}X^{n+1}/{\partial \xi})$. This can
be often done analytically. Equation (\ref{simpler})  relates $X^{n+1}$ to $\xi$ implicitly. 
We have values of $\xi$
and the corresponding values of $X^{n+1}$; to find $J$ there is no need to solve again for $X^{n+1}$;
an implicit differentiation is all that is needed. Alternately, $J$ can be found numerically,
by taking nearby values of $\xi$, redoing the iteration (which should converge in one step, because one can start
from the known value of $X^{n+1}$), and differencing. 

The expression on the
right-hand side of equation (\ref{main}) is proportional to $P(b^{n+1}|X^{n+1})P(X^{n+1}|X^n)$, with a 
proportionality constant independent of $X^n$.
When $X^{n+1}$ is sampled as just described, each value of $X^{n+1}=X^{n+1}(\xi)$ appears
with probability $\frac{1}{(2\pi)^{m/2}}\exp(-\xi^*\xi/2)/|J|$, and then the value of this expression is 
$\exp(-\xi^*\xi/2)\exp(-\Phi^{n+1})$. To get the right value of the expression on the average,
one has to give each proposed  $X^{n+1}$ the sampling weight $W=\frac{1}{(2\pi)^{m/2}}\exp(-\Phi^{n+1})|J|$,
(with another factor $P(X^n)$ if such factors are not all equal). Since $\frac{1}{(2\pi)^{m/2}}$ is a constant and the
same to every particle, we will drop it from now on.
Here we see an advantage of starting from a prechosen reference
variable $\xi$: the factor $\exp(-\xi^*\xi/2)$, which varies from sample to
sample, has been discounted in advance and does not contribute to the
non-uniformity of the weights. We shall see that the other factors can be 
expected to vary little. 

Do this for all the particles and obtain new positions with weights $W_i=\exp(-\Phi^{n+1}_i)|J_i|$, where $ \Phi^{n+1}_i, J_i$ are the values of these quantities for the $i$-th particle.  One can get rid of the weights after the fact by resampling,
i.e., for each of $M$ random numbers
$\theta_k,k=1,\dots,M$ drawn from the uniform distribution on $[0,1]$,
choose a new ${\widehat X}^{n+1}_k=X^{n+1}_i$
such that
$A^{-1}\sum_{j=1}^{i-1}W_j <  \theta_k \le A^{-1}\sum_{j=1}^i W_j$
(where $A=\sum_{j=1}^M W_j$),
and then suppress the hat.

Note also that the resampling does not have to be done at every step- for example, one can
add up the phases for a given particle and resample only when the ratio of the
largest cumulative weight $\exp(-\sum (\phi_i-\log |J_i|))$  to the smallest such weight exceeds some limit $L$ (the summation is over the  weights
accrued to a particular particle $i$ since the last resampling).
If one is worried by too many particles being close to
each other ("depletion" in the usual 
Bayesian terminology), one can divide the set of
particles into subsets of small size
and resample only inside those subsets, creating a
greater diversity. As will be seen
in the numerical results section, none of these strategies is
used here and we resample fully at every step.

The computational complexity of this construction depends on the sparseness of the matrix $\Sigma_j$, which
depends on the sparseness of  $H_{j}$ in the expression (\ref{H12}), which depends on the structure of
the function $h$ in equation (\ref{eq:observe}). In the frequently encountered situation where $h$ is diagonal, 
in the sense that each quantity
measured is a function of a single component of the vector whose dynamics are given by
equation (\ref{eq:datass}), one finds that $\Sigma_j$ and $H_{j}$ are diagonal, and the computations, including the computation
of the Jacobian $J$, are easy, whether $h$ is linear or not. The more arguments in each of the
components of the function $h$, the more labor is required.

If both equations (\ref{eq:datass}) and (\ref{eq:observe}) are linear and the
initial data are Gaussian, then the pdfs $P_n$ are Gaussian. 
We only need to find the mean and the variance of the pdf, which can be found as above by considering a single particle;
the iterations converge in one step. The resulting means and
variances are identical to those produced by the Kalman filter. If one had needed multiple particles, their weights would have been all equal. If equation (\ref{eq:datass}) is nonlinear
but equation (\ref{eq:observe}) is linear (or can be well approximated by a linear function
in each interval $(n\delta,({n+1})\delta)$), then the $P_{n+1}$ are in general not Gaussian and one
needs multiple particles. The iterations still converge in one step, and what one obtains is a
version of the forward step in a filter with an optimal importance function (as described e.g in \cite{CT1}).

The convergence of the iteration will be very briefly discussed further below. 
We have chosen the variables $\xi$ to be independent $N(0,1)$ variables,
but there is nothing sacred about this choice. The goal is to pick samples whose
probability is high, and in some contexts other choices may be better. We will discuss
those other choices when they are made in further work.

\section{Backward sampling}

In the previous section we described how to sample the pdf at time $(n+1)\delta$ given the pdf
at time $n\delta$. In general, this is not sufficient. Every observation provides information
not only about the future but also about the past- it may, for example,
tag as improbable earlier states that had seemed probable before the observation was made.
Furthermore, in non-Gaussian settings, the pdf one obtains by going directly from time $(n-1)\delta$ to step $(n+1)\delta$
by a step of duration $2\delta$ may be different from the pdf one obtains after two steps
that include an intermediate step. After one has sampled at time $(n+1)\delta$, one has
to go back, correct the past, and resample 
(this backward sampling is often misleadingly explained in the literature solely by the need to create greater diversity
among the particles).
We resample by interpolation, which we present explicitly for 
one backward step. It is quite obvious one can do that for as many backward steps
as are needed. 

Given a set of particles at time $(n+1)\delta$, after a forward step and maybe a subsequent resampling,
one can figure out where each particle $i$  was in the previous two
steps, and have a partial history for each particle $i$:
$X_i^{n-1},X_i^{n},X_i^{n+1}$ (if resamplings had
occurred,
some parts of that history may be shared among several current
particles). Knowing the first and the last members of this sequence, we 
recompute $X^n$ by interpolation, 
thus projecting information
backward one step.

The probability of the $X^{\text {new}}$ that will replace $X^n$
is the product of the three probabilities (properly normalized): the probability of the new leg
from $X^{n-1}$ to $X^n$, the probability of the resulting leg from $X^n$ to $X^{n+1}$ (the end result
being known), and the probability of the resulting observation at time $n\delta$, i.e.:
\begin{flalign}
&\exp\left(-\left(X^{\text {new}}-X^{n-1}-F_{n-1}\right)^*(G_{n-1}^*G_{n-1})^{-1}
\left(X^{\text {new}}-X^{n-1}-F_{n-1}\right)/2\right.\nonumber\\
&\left.-\left(X^{n+1}-X^{\text {new}}-F_{n}\right)^*(G_n^*G_{n})^{-1}
\left(X^{n+1}-X^n-F_{n}\right)/2
-\left(h(X^{\text{new}})-b^{n}\right)^*(Q^*Q)^{-1}
\left(h(X^{\text {new}})-b^{n}\right)/2\right).\nonumber\\
\label{new}
\end{flalign}
Here we recall that $F_{n-1}=F(X^{n-1},t^{n-1})\delta$ and  $G_{n-1}=\sqrt{\delta}G(X^{n-1},t^{n-1})$ are known from the approximation of the SDE,  $F_n$ and $G_n$ are functions of $X^{\text{new}}$, and the subscript $i$ referring to the particle has been omitted. 
This expression differs from equation (\ref{Pn+1}) by having an additional
exponential factor. 

Once again, we set up an iteration, with iterates $X_j$, that converges to $X^{\text {new}}$, and start with $X_0=0$.
We expand $h(X_{j+1})$ in a Taylor series around $X_j$, so that the last factor in the
expression (\ref{new}) becomes a quadratic in $X_{j+1}$. We complete squares so that the argument of the
exponential in (\ref{new}) can be written as $(X_{j+1}-\bar{m}_j)\Sigma_j^{-1}((X_{j+1}-\bar{m}_j)/2+\Phi_j$;
equate $(X_{j+1}-\bar{m}_j)\Sigma_j^{-1}((X_{j+1}-\bar{m}_j)/2$ to $\xi^*\xi/2$, solve to get $X_{j+1}$ as a function
of $\xi$, calculate the Jacobian, and find the weight. We do this for all the particles, and resample
as needed. This concludes the backward sampling step. Note that as a result of the backward step
and the subsequent forward step, $P_{n+1}$ depends, not only on the positions of the particles
at time $n\delta$, but also on the earlier history of the system.

\section{Sparse observations}
Consider now a situation where  we do not have observations at every time step. First, assume 
that one has observation at time $(n+1)\delta$ but not at time $n\delta$. We try to sample $X^n$ and $X^{n+1}$ 
together given the observation information at time step $(n+1)\delta$. 
Consider the $i$-th particle. 
Suppose we are given the  vector $X^{n-1}_i$ for that particle.
Suppress again the 
particle
index $i$. 
The joint probability density $P_{n,n+1}$ of $X^n$ and $X^{n+1}$ given $X^{n-1}$ is
\begin{flalign}
&P_{n,n+1}(X^n,X^{n+1})\nonumber\\
=&Z^{-1}\exp\left(-\left(X^{n}-X^{n-1}-F_{n-1}\right)^*(G_{n-1}^*G_{n-1})^{-1}\left(X^{n}-X^{n-1}-F_{n-1}\right)/2\right.\nonumber\\
&\left.-\left(X^{n+1}-X^n-F_n\right)^*(G_n^*G_n)^{-1}\left(X^{n+1}-X^n-F_n\right)/2
-\left(h(X^{n+1})-b^{n+1}\right)^*(Q^*Q)^{-1}\left(h(X^{n+1})-b^{n+1}\right)/2\right),\nonumber\\
\label{bPn+1}
\end{flalign}
where 
$Z$ is the normalization constant.  We recall that $F_{n-1}=F(X^{n-1},t^{n-1})\delta$, $G_{n-1}=\sqrt{\delta}G(X^{n-1},t^{n-1})$ are known from the approximation of the SDE,  $F_n$ and $G_n$ depend on $X^{n}$.

In the now familiar sequence of steps, 
we pick two independent samples $\xi_n$ and $\xi_{n+1}$, each with probability density
$\exp(-\xi^*\xi/2)/(2\pi)^{m/2}$, and try to solve the equation
\begin{flalign}
&{\xi_n}^*\xi_n/2+\xi_{n+1}^*\xi_{n+1}/2\nonumber\\
=&\left(X^{n}-X^{n-1}-F_{n-1}\right)^*(G_{n-1}^*G_{n-1})^{-1}\left(X^{n}-X^{n-1}-F_{n-1}\right)/2\nonumber\\
&+\left(X^{n+1}-X^n-F_n\right)^*(G_n^*G_n)^{-1}\left(X^{n+1}-X^n-F_n\right)/2
+\left(h(X^{n+1})-b^{n+1}\right)^*(Q^*Q)^{-1}\left(h(X^{n+1})-b^{n+1}\right)/2,
\label{bargs}
\end{flalign}
to obtain $X^n$ and $X^{n+1}$ as  functions of $\xi_n$ and $\xi_{n+1}$.

We define 
a sequence of approximations $X^n_j$ and $X^{n+1}_j$  which converge to $X^n$ and $X^{n+1}$, respectively;
set $X^n_0=0$ and $X_0^{n+1}=0$, and at each iteration find $X^n_{j+1}$ and $X^{n+1}_{j+1}$ given $X^n_j$ and $X^{n+1}_j$.
First, expand the function $h$ in the observation equation (\ref{eq:observe}) in Taylor series around $X^{n+1}_j$:
\begin{equation}\label{bobs}
h(X^{n+1}_{j+1})=h(X^{n+1}_j)+H^{n+1}_j\cdot(X^{n+1}_{j+1}-X^{n+1}_j),
\end{equation}
where $H^{n+1}_j$ is a Jacobian matrix evaluated at $X^{n+1}_j$.
The observation equation (\ref{eq:observe}) is approximated as:
\begin{equation}
z_j^{n+1}=H^{n+1}_jX^{n+1}_{j+1}+QW^{n+1},
\label{blincon}
\end{equation}
where $z^{n+1}_j=b^{n+1}-h(X^{n+1}_j)+H^{n+1}_jX^{n+1}_j$.

Let $F_{n,j}=F(X^n_j,t^n)\delta$ and $G_{n,j}=\sqrt{\delta}G(X^n_j,t^n)$. 
The right side of equation (\ref{bargs}) can be approximated as:
\begin{flalign}
&\left(X_{j+1}^{n}-X^{n-1}-F_{n-1}\right)^*(G_{n-1}^*G_{n-1})^{-1}\left(X_{j+1}^{n}-X^{n-1}-F_{n-1}\right)/2\nonumber\\
+&\left(X_{j+1}^{n+1}-X_{j+1}^n-F_{n,j}\right)^*(G_{n,j}^*G_{n,j})^{-1}\left(X_{j+1}^{n+1}-X_{j+1}^n-F_{n,j}\right)/2\nonumber\\
+&\left(H_j^{n+1}X_{j+1}^{n+1}-z_j^{n+1}\right)^*(Q^*Q)^{-1}\left(H_j^{n+1}X_{j+1}^{n+1}-z_j ^{n+1}\right)/2.\nonumber\\
\label{bH121}
\end{flalign}
We first combine the last two terms in  (\ref{bH121}) and obtain 
\begin{flalign}
&\left(X_{j+1}^{n+1}-X_{j+1}^n-F_{n,j}\right)^*(G_{n,j}^*G_{n,j})^{-1}\left(X_{j+1}^{n+1}-X_{j+1}^n-F_{n,j}\right)/2
+\left(H_{n+1}X_{j+1}^{n+1}-z^{n+1}\right)^*(Q^*Q)^{-1}\left(H_{n+1}X_{j+1}^{n+1}-z^{n+1}\right)/2\nonumber\\
=&\left(X^{n+1}_{j+1}-\bar{m}_j^{n+1}\right)^*(\Sigma_j^{n+1})^{-1}\left(X^{n+1}_{j+1}-\bar{m}_j^{n+1}\right)/2+\Phi^{n+1}_j,
\label{bH122}
\end{flalign}
where 
$$
(\Sigma_j^{n+1})^{-1}=(G_{n,j}^*G_{n,j})^{-1}+(H_j^{n+1})^*(Q^*Q)^{-1}H_j^{n+1},
$$
$$ \bar{m}_j^{n+1}=\Sigma_j^{n+1}\left((G_{n,j}^*G_{n,j})^{-1}(X^n_{j+1}+F_{n,j})+(H_j^{n+1})^*(Q^*Q)^{-1}z_j^{n+1}\right),
$$
$$
K_j^{n+1}=H_j^{n+1}G_{n,j}^*G_{n,j}(H_j^{n+1})^*+Q^*Q,
$$
and
$$\Phi^{n+1}_j=\left(z_j^{n+1}-H_j^{n+1}(X^n_{j+1}+F_{n,j})\right)^*(K_j^{n+1})^{-1}\left(z_j^{n+1}-H_j^{n+1}(X^n_{j+1}+F_{n,j})\right)/2.
$$
We combine the first term in (\ref{bH121}) and the second term in (\ref{bH122}) and obtain
\begin{eqnarray}
&&\left(X_{j+1}^{n}-X^{n-1}-F_{n-1}\right)^*(G_{n-1}^*G_{n-1})^{-1}\left(X_{j+1}^{n}-X^{n-1}-F_{n-1}\right)/2+\Phi^{n+1}_j\nonumber\\
&=&\left(X_{j+1}^{n}-X^{n-1}-F_{n-1}\right)^*(G_{n-1}^*G_{n-1})^{-1}\left(X_{j+1}^{n}-X^{n-1}-F_{n-1}\right)/2\nonumber\\
&&+\left(z_j^{n+1}-H_j^{n+1}(X^n_{j+1}+F_{n,j})\right)^*(K_j^{n+1})^{-1}\left(z_j^{n+1}-H_j^{n+1}(X^n_{j+1}+F_{n,j})\right)/2\nonumber\\
&=&\left(X^n_{j+1}-\bar{m}_j^n\right)^*(\Sigma_j^n)^{-1}\left(X^n_{j+1}-\bar{m}_j^n\right)/2+\Phi^n_j,
\label{bH123}
\end{eqnarray}
where
$$
(\Sigma_j^n)^{-1}=(G_{n-1}^*G_{n-1})^{-1}+(H_j^{n+1})^*(K_jj^{n+1})^{-1}H_j^{n+1},
$$
$$ \bar{m}_j^n=\Sigma_j^n\left((G_{n-1}^*G_{n-1})^{-1}(X^{n-1}+F_{n-1})+(H_j^{n+1})^*(K_j^{n+1})^{-1}(z_j^{n+1}-H_j^{n+1}F_{n,j})\right),
$$
$$
K_j^n=H_j^{n+1}G_{n-1}^*G_{n-1}(H_j^{n+1})^*+K_j^{n+1},
$$
and
\begin{eqnarray*}
\Phi_j^n=\left(z_j^{n+1}-H_j^{n+1}(F_{n,j}+X^{n-1}+F_{n-1})\right)^*(K_j^n)^{-1}
\left(z_j^{n+1}-H_j^{n+1}(F_{n,j}+X^{n-1}+F_{n-1})\right)/2.
\end{eqnarray*}
Combining (\ref{bargs}), (\ref{bobs}), (\ref{bH121}), (\ref{bH122}), and (\ref{bH123}), we try to solve
\begin{flalign}
&{\xi_n}^*\xi_n/2+\xi_{n+1}^*\xi_{n+1}/2\nonumber\\
=&\left(X^{n+1}_{j+1}-\bar{m}_j^{n+1}\right)^*(\Sigma_j^{n+1})^{-1}\left(X^{n+1}_{j+1}-\bar{m}_j^{n+1}\right)/2+\left(X^n_{j+1}-\bar{m}_j^n\right)^*(\Sigma_j^n)^{-1}\left(X^n_{j+1}-\bar{m}_j^n\right)/2+\Phi^n_j.\label{bb}
\end{flalign}
We now solve for $X^n_{j+1}$ and $X^{n+1}_{j+1}$ as  functions of $\xi_n$ and $\xi_{n+1}$, ignoring the remainders $\Phi^n_j$, i.e. we solve
the simpler equations
\begin{equation}
(X^k_{j+1}-\bar{m}_j^k)^*(\Sigma_j^k)^{-1}(X^k_{j+1}-\bar{m}_j^k)/2=\xi_k^*\xi_k/2, \quad k=n,n+1
\label{bsimpler}
\end{equation}
If the sequences $X^n_j$ and $X^{n+1}_j$ converge to limits, call the limits $X^n$ and  $X^{n+1}$.  
In the limit, the approximate equation (\ref{blincon})
converges to the full observation equation (\ref{eq:observe}). The remainders
$\Phi^n_j$ and $\Phi^{n+1}_j$ also converge to limits $\Phi^n$ and $\Phi^{n+1}$. Equation
(\ref{bargs}) becomes:
\begin{eqnarray}
&&\xi_n^*\xi_n/2+\xi_{n+1}^*\xi_{n+1}/2+\Phi^n\nonumber\\
&=&\left(X^{n}-X^{n-1}-F_{n-1}\right)^*(G_{n-1}^*G_{n-1})^{-1}\left(X^{n}-X^{n-1}-F_{n-1}\right)/2\nonumber\\
&&+\left(X^{n+1}-X^n-F_n\right)^*(G_n^*G_n)^{-1}\left(X^{n+1}-X^n-F_n\right)/2
+(h(X^{n+1})-b^{n+1})(Q^*Q)^{-1}(h(X^{n+1})-b^{n+1})/2.\nonumber\\
\label{bwemiss}
\end{eqnarray}
Multiply by $-1$ and exponentiate:
\begin{eqnarray}
&&\exp(-\xi_n^*\xi_n/2)\exp(-\xi_{n+1}^*\xi_{n+1}/2)\exp(-\Phi^n)\nonumber\\
&=&\exp\left(\left(X^{n}-X^{n-1}-F_{n-1}\right)^*(G_{n-1}^*G_{n-1})^{-1}\left(X^{n}-X^{n-1}-F_{n-1}\right)/2\right.\nonumber\\
&+&\left.\left(X^{n+1}-X^n-F_n\right)^*(G_n^*G_n)^{-1}\left(X^{n+1}-X^n-F_n\right)/2
+\left(h(X^{n+1})-b^{n+1})^*(Q^*Q)^{-1}(h(X^{n+1})-b^{n+1}\right)/2\right).\nonumber\\
\label{bmain}
\end{eqnarray}

As before, one has to give each proposed $X^n$ and $X^{n+1}$ the sampling 
weight $W=\exp(-\Phi^n)|J|$, 
where $J$ is the Jacobian 
$J=\det (\partial (X^n,X^{n+1})/\partial (\xi_n,\xi_{n+1}))$ which must be computed.  
One does this for all
particles and resamples as needed. This process can be generalized if one
wishes to sample at more times between observations. One  should also note that
the procedure just described may make the evaluation of Jacobians significantly more onerous, but still often tractable.

The construction of this paragraph is important because many data sets one
tries to assimilate are indeed sparse, and also for the following reason.
We have not provided in this present paper a discussion of the convergence of the iterations
we use. This convergence depends on the structure of
the underlying SDE, on the scheme used to approximate it, and on the specific ways one solves for the new increments in terms of the reference variables $\xi$,  and cannot
be analyzed without considering these specifics. In our previous paper \cite{CT1}
we analyzed a special case and found that there the convergence depended on the
size of the time step. We conjecture that this happens frequently.
The present section provides a way to decrease the time step as a device
for repairing diverging iterations without much additional thought.  

\section{Example 1}
We apply our filter to a prototypical marine ecosystem model studied 
in \cite{Dow1}. We set the main parameters equal to the ones in \cite{Dow1}; however, we will also present some results with a range of noise variances to make
a particular point. We did the data assimilation with the filter described above, without back sampling, and also by the a standard particle filter SIR (Sampling importance resampling), see \cite{Ar2}.

The model involves
four state variables: phytoplankton P (microscopic plants), 
zooplankton Z (microscopic animals), nutrients N (dissolved inorganics), and 
detritus D (particulate organic non-living matter). At the initial time $t=0$ we have $P(0)=0.125$, $Z(0)=0.00708$, $N(0)=0.764$, and $D(0)=0.136$. The system is described by the following nonlinear ordinary differential equations, explained in \cite{Dow1}:
\begin{eqnarray}
\frac{dP}{dt}&=&\frac{N}{0.2+N}\gamma P-0.1P-0.6\frac{P}{0.1+P}Z+N(0,\sigma^2_P)\nonumber\\
\frac{dZ}{dt}&=&0.18\frac{P}{0.1+P}Z-0.1Z+N(0,\sigma^2_Z)\nonumber\\
\frac{dN}{dt}&=&0.1D+0.24\frac{P}{0.1+P}Z-\gamma P\frac{N}{0.2+N}+0.05Z+N(0,\sigma^2_N)\nonumber\\
\frac{dD}{dt}&=&-0.1D+0.1P+0.18\frac{P}{0.1+P}Z+0.05Z+N(0,\sigma^2_D),
\label{ODE}
\end{eqnarray}
where the parameter $\gamma$, the `` growth rate",  is determined by the equations given by
$$\gamma_t = 0.14 + 3\Delta \gamma_t,\quad \Delta\gamma_t =0.9\Delta \gamma_{t-1}+N(0,\sigma^2_\gamma).$$
The variances of the noise terms are: $\sigma_P^2=(0.01P(0))^2$, $\sigma_Z^2=(0.01Z(0))^2$,  $\sigma_N^2=(0.01N(0))^2$, $\sigma_D^2=(0.01D(0))^2$, and  $\sigma_\gamma^2=(0.01)^2$.

The observations were obtained from NASA's 
SeaWiFS satellite ocean color images. These observations provide a time series for 
phytoplankton;  the relation between the observations $P(t)_{\mbox{obs}}$ (corresponding to the vector $b^n$ in the earlier discussion) and the solution $P(t)$ of the equation of the first equation in (\ref{ODE}) is assumed to be: 
$$\log P(t)_{\mbox{obs}}=\log P(t) +N(0,\sigma^2_{\mbox{obs}}),$$
where $\sigma^2_{\mbox{obs}}=0.3^2$. Note that this observation equation is not linear. There are 190 data points distributed
 from late 1997 to mid 2002. The sample intervals ranged
 from a week to a 
month or more, for details see \cite{Dow1}. 
As in \cite{Dow1}, we discretize the system (\ref{ODE}) by an Euler method
with $\Delta t=1$ day and prohibit the state variables from dropping below 1 percent of their initial values. 

We have compared our filter and SIR  in three sets of numerical experiments, all with the same initial values as listed above. In each case we attempted
to find a trajectory of the system consistent with the fixed data, and observed
how well we succeeded. 
In the first set of the experiments, we used 100 particles and take $\sigma_P^2=(0.01P(0))^2$ as in \cite{Dow1}. In this case, the (assumed) variance
of the system is much smaller than the (assumed) variance of the observations; the particle paths are bunched close together, and the results from our filter and from 
SIR are quite close, see Figure 1, where we plotted the $P$ component of
the reconstructed solution as well as the corresponding data. 

In the second set of the experiments, we still used 100 particle but assumed $\sigma_p^2=(P(0))^2$. The variance of the system is now comparable to
the variance of the observation. For SIR, after resampling, the number of the distinct particles is smaller than in the first case, as a result of the loss of diversity after resampling when the weights are very different from each other, see Table 1, where we exhibit the average number of distinct particles left after each resample; there is a resample after each step.
Remember that there is some loss of diversity in resampling even if all the
weights are equal.
With 100 
particles, the filtered results with SIR are still comparable to those with our filter. See Figure 2.

In the third set of the experiments, we used only 10 particles and kept $\sigma_p^2=(P(0))^2$.  As one could have foreseen, our filter does better than
SIR, see Figure 3. One should remember however that we are working with a low
dimensional problem where the differences between filters are not expected
to be very significant; the cost if 100 particles is not prohibitive.

\begin{table}
\caption{The number of distinct particles after resampling with different system variances and different numbers of particles }
\center{
\begin{tabular}{ |c|c|c|c| }\hline\label{T1}
 $\sigma_p$ & {\#} particle  &\multicolumn{2}{c|}{average {\#} particles left after resampling} \\\hline
 &   &   SIR &Our filter \\\hline
$0.01P(0)$ &100&61& 61  \\\hline
$P(0)$     &100&19& 63   \\\hline
$P(0)$     &10&2.2& 6.3    \\\hline
\end{tabular}
}
\end{table}

\begin{center}
\begin{figure}
\caption{Results with $\sigma_P^2=(0.01P(0))^2$ and 100 particles}
\includegraphics[scale=0.5,angle=0]{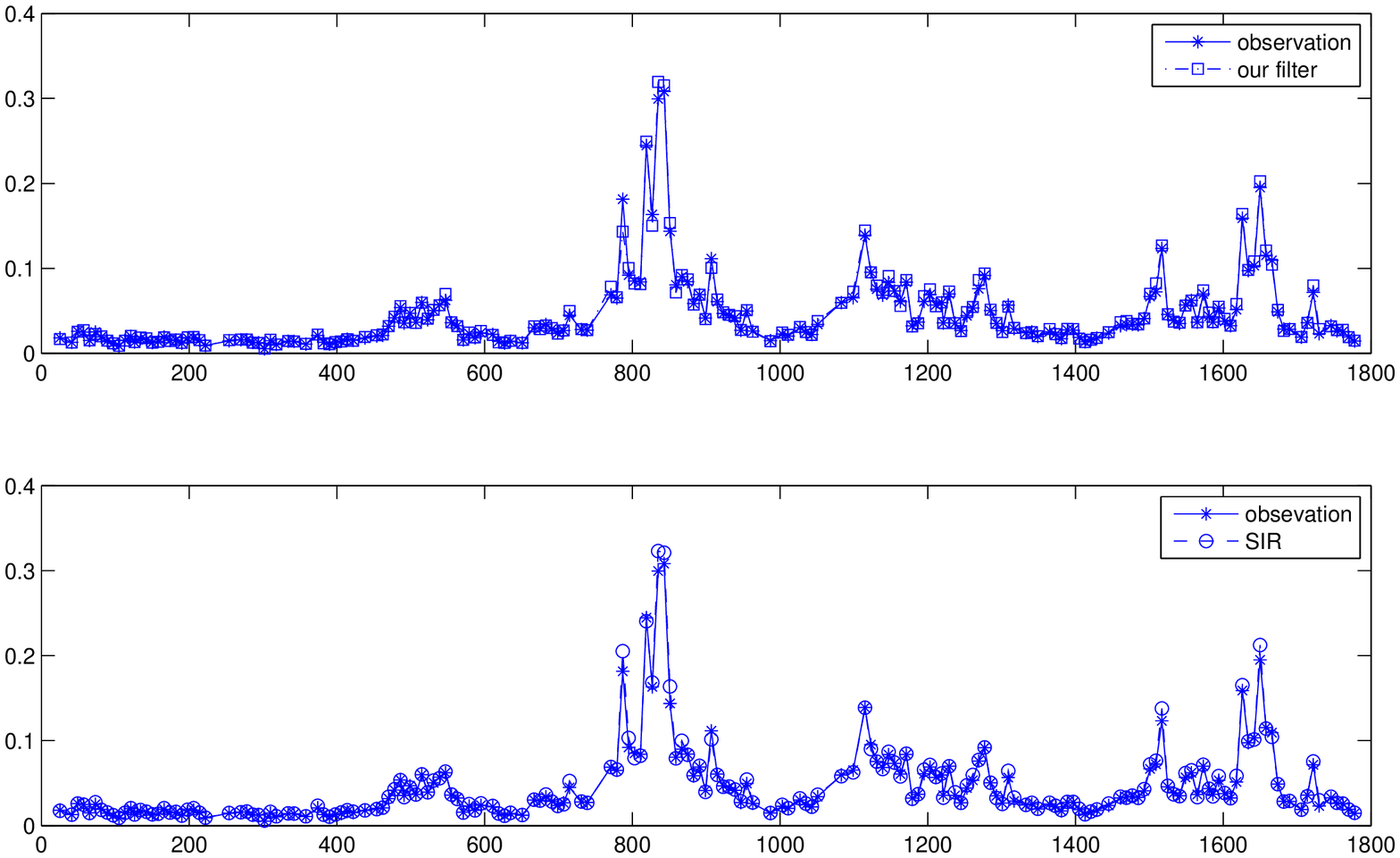}
\end{figure}
\begin{figure}
\caption{Results with $\sigma_P^2=P(0)^2$ and 100 particles}
\includegraphics[scale=0.5,angle=0]{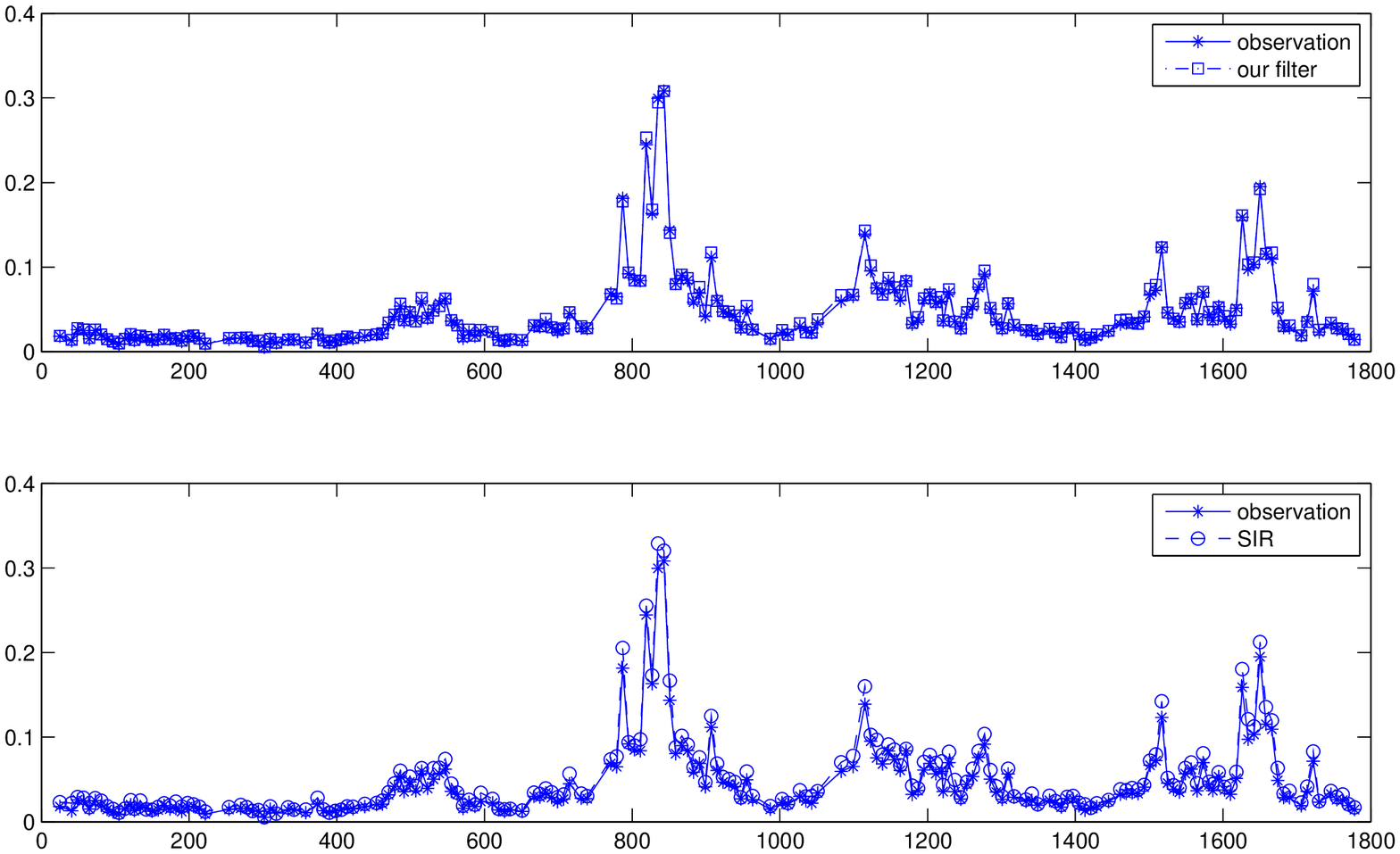}
\end{figure}
\begin{figure}
\caption{Results with $\sigma_P^2=P(0)^2$ and 10 particles}
\includegraphics[scale=0.9,angle=90]{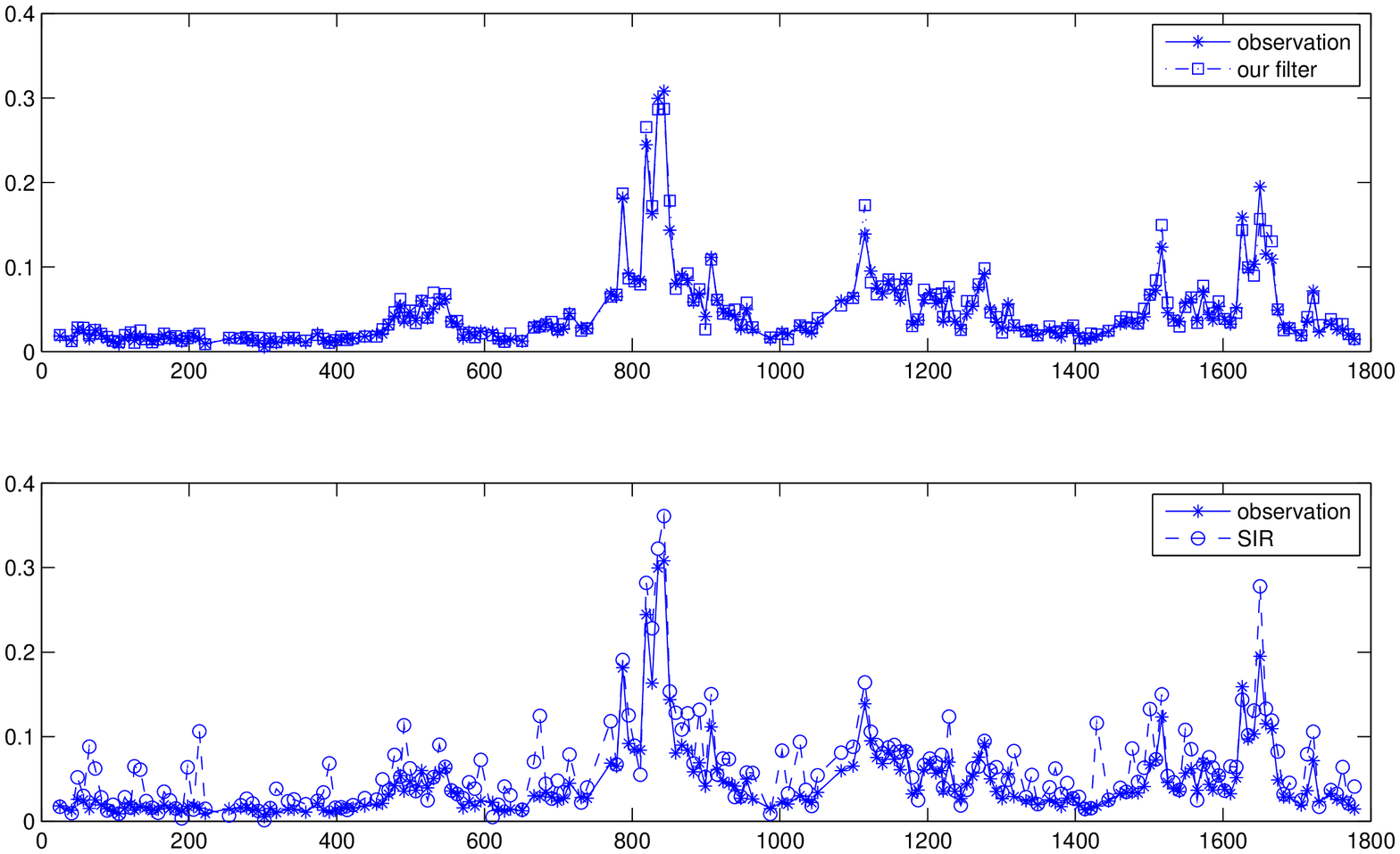}
\end{figure}
\end{center}

\section{Example 2}

We consider next a simple high dimensional example, used in \cite{Sny} to show 
how particle filters fail when the number of dimensions is large. 
We assume that each component of $X^n$ is an independent Gaussian 
with zero 
mean and unit variance. This is equivalent to 
 taking $\delta=1$, $F(X^n,\delta)=0$, $G(X^n,t^n)=I$ in equation (\ref{appeq}),
 and eliminating the $X^n$ term. We have
\[X^n=V^n.\]
Each component of $X^n$ is observed individually, so that 
\[b^n=X^n+W^n.\]
We implement our filter with these particular choices. At the $j$-th iteration, $H_j=I$ in equation (\ref{lineeqn}) and 
$z_j=b^{n+1}$ in equation (\ref{lincon}). Therefore, we have
$\Sigma_j^{-1}=2I$, $\bar{m}_j=b^{n+1}/2$, and $\Phi_j=(b^{n+1})^*b^{n+1}/4$, in equation (\ref{H12}). The iterations
converge in one step and all the particles have the same weights.

However, with SIR the weights are uneven. 
We ran the SIR filter 1000 times, with a 1000 particles each time; in each run we
normalized the weights so that add up to one, and we  recorded the maximum weight. 
In Figures 4 we display a histogram of  these recorded maximum weights.
As one can observe, when
the number of dimensions is large, most of time, a single particle in each run hogs
all the probability, and this version of SIR fails.
\begin{center}
\begin{figure}
\caption{Histogram of the SIR normalized maximum particle weights with 1000 runs for $100$ dimensions}
\includegraphics[scale=0.6]{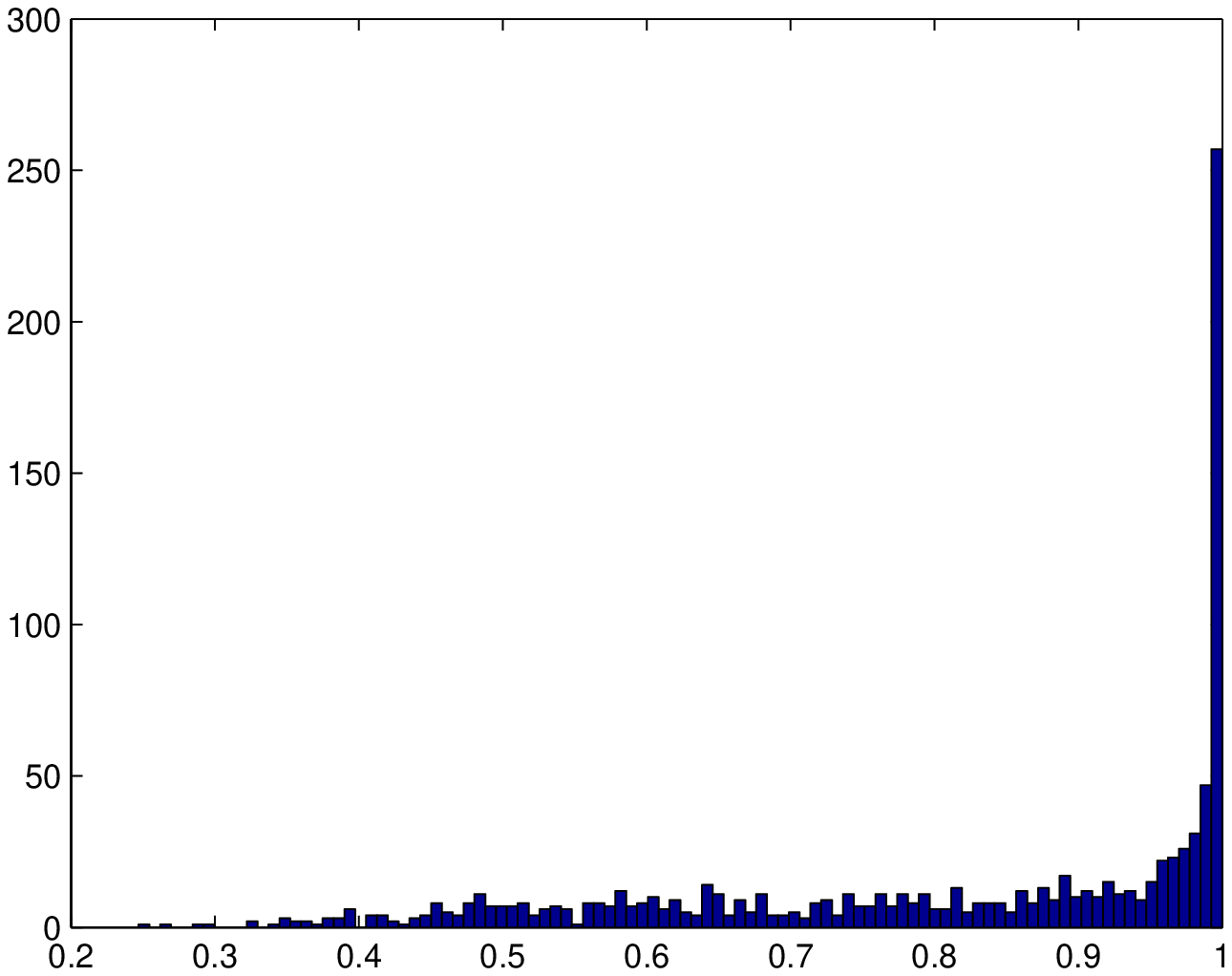}
\end{figure}
\end{center}

\section{Conclusions}
We have presented a general form of the iteration and interpolation process used in our new implicit nonlinear particle
filter. 
The goal is to aim particle paths sharply so that fewer are needed.
We conjecture that there is no general way to reduce the variability of the weights in particle sampling further than we have.
We also presented additional simple examples that illustrate the potential
of this new sampling. These examples are simple in that one is low-dimensional, while
the second is linear so that other effective ways of sampling it do
exist. High-dimensional nonlinear problems where our filter may be
indispensable will be presented elsewhere, in the context of specific applications.

\section{Acknowledgments}
We would like to thank Prof. J. Goodman, who urged us to write a more
general version of our previous work and suggested some notations and nomenclature,
 Prof. R. Miller, who suggested that we try Dowd's model plankton problem as a first
step toward an ambitious joint effort and helped us set it up, and Prof. M. Dowd, who kindly made the data available.  
This work was supported in part by the Director, Office of Science,
Computational and Technology Research, U.S.\ Department of Energy under
Contract No.\ DE-AC02-05CH11231, and by the National Science Foundation under grant
 DMS-0705910.

\bibliography
{filter}
\bibliographystyle{plain}

\end{document}